\documentclass[journal]{IEEEtran}

\usepackage{float}
\usepackage{caption}
\usepackage{subcaption}
\usepackage{tikz,pgfplots,pgfplotstable,booktabs}
\usepackage[official]{eurosym}
\usepackage{amsmath}
\usepackage{amssymb}
\allowdisplaybreaks[1]

\begin{document}

\title{Is learning for the unit commitment problem \\ a low-hanging fruit?}
\author{S. Pineda and J. M. Morales
\thanks{S. Pineda is with the Dep.
of Electrical Engineering, Univ. of Malaga, Spain. E-mail: spinedamorente@gmail.com. J. M. Morales is with the Dep. of Applied Mathematics, Univ. of Malaga, Spain. E-mail: juan.morales@uma.es.}
\thanks{This work was supported in part by the Spanish Ministry of Science and Innovation through project PID2020-115460GB-I00, by the Andalusian Regional Government through project P20-00153, and by the Research Program for Young Talented Researchers of the University of Málaga under Project B1-2019-11. This project has also received funding from the European Research Council (ERC) under the European Union’s Horizon 2020 research and innovation programme (grant agreement No 755705).}}

\maketitle

\begin{abstract}
The blast wave of machine learning and artificial intelligence has also reached the power systems community, and amid the frenzy of methods and black-box tools that have been left in its wake, it is sometimes difficult to perceive a glimmer of Occam's razor principle. In this letter, we use the unit commitment problem (UCP), an NP-hard mathematical program that is fundamental to power system operations, to show that \emph{simplicity} must guide any strategy to solve it, in particular those that are based on \emph{learning} from past UCP instances. To this end, we apply a naive algorithm to produce candidate solutions to the UCP and show, using a variety of realistically sized power systems, that we are able to find optimal or quasi-optimal solutions with remarkable speedups. Our claim is thus that any sophistication of the learning method must be backed up with a statistically significant improvement of the results in this letter.
\end{abstract}


\IEEEpeerreviewmaketitle

\section{Introduction}\label{sec:intro}

\IEEEPARstart{T}{he} unit commitment problem (UCP) is currently one of the most fundamental mathematical tools to operate power systems. The UCP determines the on/off commitment status and power dispatch of generating units to satisfy electricity demand at a minimum cost while complying with the technical limits of generation and transmission assets \cite{sen1998optimal}. The UCP is usually formulated as a mixed-integer optimization problem that is proven to be NP-hard \cite{bendotti2019complexity} and therefore, the technical literature includes several methods to trim down its computational burden \cite{saravanan2013solution, chen2016improving}. Existing strategies comprise formulation tightening \cite{pandzic2013comparison,tejada2020which}, decomposition techniques \cite{fu2005security} and constraint screening \cite{zhai2010fast}. These methods, however, overlook the fact that slight variations of the same UCP  are to be solved everyday and therefore, learning from the past may be a powerful weapon to tackle future UCP instances. 

In the same vein, some learning-based methodologies have been recently proposed to reduce the computational burden of the UCP using information about previously solved instances.  References \cite{yang2021machine,ruan2021review} review  current learning-based methods for power system problems.
In particular, the authors of \cite{dalal2018unit} learn a proxy of the UCP's solution to be used in long-term planning assessments. Despite being simple and fast (with speedups up to 260x), the proposed strategy involves an average optimality error of 3.5\% for the IEEE-96RTS test system, which precludes its use in short-term operation. Using Classification and Regression Trees (CART) or Random Forest (RF), reference \cite{lin2019approximate} presents a methodology to find the relationships between the solutions of the original and relaxed versions of the UCP and reports average optimality errors of between 0.14\% and 0.23\%. A supervised classification procedure is proposed in \cite{pineda2020data} to learn the transmission capacity constraints that can be screened out of the UCP. Using a 2000-bus system, the authors report a speedup of 19x if retrieving the original solutions must be guaranteed, and a speedup of 43x for an average suboptimality error of 0.04\%. The authors of \cite{yang2020integrated} describe a sophisticated methodology to cluster decision variables depending on the difficulty of the UCP instance to be solved. They achieve speedups of between 1.5x and 2x, an average optimality error of 0.1\%, and a maximum optimality error of 1\% for the IEEE 118-bus system. The authors of \cite{chen2020distributed} use the UCP solution of the previous day UC as a warm start strategy and report speedups of 2x for the MISO power system. This speedup can be increased to 2x-12x if additional historical information is considered, as shown in \cite{chen2021high}. In \cite{mohammadi2021machine}, the authors use machine learning to determine unnecessary constraints of the stochastic unit commitment and attain speedups of 14x for a 500-bus power system. Finally, reference \cite{xavier2021learning} proposes different strategies to learn initial solutions with speedups of 4.3x while retrieving the optimal solution of the UCP. The authors also discuss a learning-based procedure to learn the relationships among binary variables, which yields speedups of 10.2x but no optimality guarantees. These methodologies are tested on a set of large-sized networks.

It is apparent that when it comes to learning for the UCP, one may easily get lost in a myriad of methods and approaches, all of which promise reasonable computational savings. The takeaway message is that, despite the high theoretical complexity of the UCP, having access to previous instances may significantly reduce the solution task in practice, since today's commitment decisions are probably very similar to those made yesterday, one week ago, or even one year ago. Given the high potential of using historical data to reduce the computational burden of the UCP, the questions that naturally arise are: Are existing learning methods an actual breakthrough in the solution of the UCP or are they just picking the low-hanging fruit? Do sophisticated learning methods to solve the UCP substantially outperform painless naive learning strategies? How should we benchmark the performance of learning methods to solve the UCP? In this letter, we aim to answer these questions by learning the UCP's solution through a straightforward K-nearest neighbor procedure. The purpose of this letter is, by no means, to propose a learning-based methodology that outperforms existing ones under certain conditions. What we claim is that the performance of existing and upcoming learning-based methods to solve the UCP should be thoroughly compared with some painless naive methods such as the one we suggest, in order to justify the hassle from the  increased complexity and sophistication, and the loss of transparency.

\section{A Naive Learning Method}\label{sec:learning}

The unit commitment problem can be generally formulated as the following optimization program:
\begin{subequations}
\begin{align}
\underset{\textbf{u},\textbf{y}}{\min} \enskip & f(\textbf{u},\textbf{y}) \\
 \text{s.t.} \enskip & g_j(\textbf{u},\textbf{y},\textbf{d}) \leq 0, \forall j 
\end{align} \label{eq:general_uc}
\end{subequations}
\noindent where $\textbf{u}$ and $\textbf{y}$ denote, respectively, the binary and continuous variables, $\textbf{d}$ represents the input parameters such as net demand throughout the network, and $f(\cdot),g_j(\cdot)$ are the objective function and the technical constraints, in that order. Under some mild assumptions, model \eqref{eq:general_uc} becomes a mixed-integer quadratic programming problem that can be solved using optimization solvers at a usually high computational cost. With some abuse of notation, we express the solution of \eqref{eq:general_uc} as a function of the input parameters, namely, $\textbf{u}^{\rm UC}(\textbf{d}),\textbf{y}^{\rm UC}(\textbf{d})$. If binary variables $\textbf{u}$ are fixed to given values $\Tilde{\textbf{u}}$, model \eqref{eq:general_uc} becomes an optimal power flow (OPF) problem, which is easier to solve and whose optimal solution is denoted as $\textbf{y}^{\rm OPF}(\Tilde{\textbf{u}},\textbf{d})$. 
    
Suppose we have access to a sample set of optimal solutions of problem \eqref{eq:general_uc} for different input parameters, referred to as $S=\{(\textbf{d}_i,\textbf{u}^{\rm UC}_i,C^{\rm UC}_i,T^{\rm UC}_i)\}_{i\in\mathcal{I}}$, where $\textbf{u}^{\rm UC}_i = \textbf{u}^{\rm UC}(\textbf{d}_i)$, and $C^{\rm UC}_i,T^{\rm UC}_i$ are, respectively, the objective function and the computational time of problem \eqref{eq:general_uc} for instance $i \in \mathcal{I}$. Intuitively, the naive learning methodology we propose as a benchmark consists in fixing the binary variables to those of close past instances to solve several OPFs in parallel and select the one with the lowest cost. Using a leave-one-out procedure, the learning strategy for each instance $i\in\mathcal{I}$ runs as follows:
\begin{itemize}
    \item[1)] Compute the distance between the input parameters $\textbf{d}_i$ and those of the remaining instances using a norm, for example, $||\textbf{d}_i-\textbf{d}_{\Tilde{i}}||_2, \forall \Tilde{i}\in\mathcal{I},\Tilde{i}\neq i$.
    \item[2)] Find the set of the $K$-nearest neighbors to  $i$ with the lowest distances computed in step 1), denoted as $\mathcal{I}^K_i$.
    \item[3)] Solve the optimal power flow for input parameters $\textbf{d}_i$ and binary variables fixed to $\textbf{u}^{\rm UC}_{\Tilde{i}}$. That is, compute $\textbf{y}^{\rm OPF}(\textbf{u}^{\rm UC}_{\Tilde{i}},\textbf{d}_i), \forall \Tilde{i} \in \mathcal{I}^K_i$ and denote the objective function and solving time as $C^{\rm OPF}_{i\Tilde{i}}$ and $T^{\rm OPF}_{i\Tilde{i}}$, respectively. 
    \item[4)] Among all feasible problems solved in step 3), approximate the cost of the UCP for instance $i$ as $\widehat{C}^{\rm UC}_i = \min_{\Tilde{i}}C^{\rm OPF}_{i\Tilde{i}}$. The lowest suboptimality gap is thus computed as $\Delta_i = (\widehat{C}^{\rm UC}_i-C^{\rm UC}_i)/C^{\rm UC}_i$.
    \item[5)] Problems in step 3) are solved in parallel and the speedup factor is thus computed as $S_i = T^{\rm UC}_i/(\max_{\Tilde{i}}(T^{\rm OPF}_{i\Tilde{i}})+T^{L}_i)$, where $T^{L}_i$ is the learning time of steps 1) and 2).
\end{itemize}

Once steps 1)-5) are run for each instance, we determine the average suboptimality gap $\overline{\Delta}$, the maximum suboptimality gap $\Delta^{\max}$, the average speedup $\overline{S}$, and the number of instances for which the $K$ problems solved in step 3) are infeasible $N^{\rm IN}$.

\section{Numerical results}

In this section, we provide the results obtained by the learning method described in Section \ref{sec:learning} for nine large-scale European test systems used in \cite{xavier2021learning} and available for download at \cite{alison2020dataset}. The numbers of buses, units and lines of each system are collated in columns 2-4 of Table \ref{tab:results}.

\begin{table*}[]
\renewcommand\arraystretch{1.3}
\centering
\begin{tabular}{ccccccccccccc}
\hline
System&Buses&Units&Lines&$\overline{\Delta}$&$\Delta^{\max}$& $< 0.01\%$ &$0.01-0.02\%$&$0.02-0.05\%$&$0.05-0.1\%$&$>0.1\%$&$N^{\rm IN}$&$\overline{S}$ \\
\hline
1888rte&1888&297&2531&0.0174&0.2394&230&131&109&20&9&1&116.5x\\
1951rte&1951&391&2596&0.0382&0.3759&47&116&217&85&27&8&150.4x\\
2848rte&2848&547&3776&0.0186&0.1332&179&138&159&14&8&2&132.6x\\
3012wp&3012&502&3572&0.0485&0.4864&37&78&212&132&36&5&188.8x\\
3375wp&3374&596&4161&0.1256&0.8073&9&14&102&164&198&13&215.9x\\
6468rte&6468&1295&9000&-0.0001&0.0175&498&2&0&0&0&0&41.2x\\
6470rte&6470&1330&9005&-0.0016&0.0187&496&4&0&0&0&0&171.9x\\
6495rte&6495&1372&9019&-0.0001&0.0481&496&3&1&0&0&0&41.0x\\
6515rte&6515&1388&9037&-0.0009&0.0133&497&3&0&0&0&0&101.7x\\ 
\hline
\end{tabular}
\caption{Numerical results}
\label{tab:results}
\end{table*}

The performance of the naive learning method is illustrated using the solution of 500 UCP instances that differ on the 24-hour load profile. According to the procedure proposed in \cite{xavier2021learning}, each 24-hour load profile is randomly generated as follows:
\begin{itemize}
    \item[1)] Let $\overline{P}_g$ denote the  capacity of each unit $g$. The peak demand $\overline{D}$ is randomly generated as $\overline{D} = 0.6 \times \left(\sum_g \overline{P}_g\right) \times \text{unif}(0.925,1.075)$, where $\text{unif}(a,b)$ represents a random sample from a uniform distribution in the interval $[a,b]$.
    \item[2)] Let $\overline{\beta}_b$ denote the nominal load distribution factor of bus $b$. The distribution factor $\beta_b$ is randomly generated as $\beta_b=\overline{\beta}_b \times \text{unif}(0.9,1.1)$, which is subsequently normalized as $\beta^{\rm N}_b = \beta_b/\sum_b \beta_b$ to satisfy that $\sum \beta^{\rm N}_b = 1$.
    \item[3)] Let $t \in \mathcal{T}$ denote the index for the time periods of the UCP, and let $\mu_t$ and $\sigma_t$ represent, respectively, the average and standard deviation of the ratio between the aggregated load level for two consecutive hours $D_{t}/D_{t-1}$. The hourly variation factors $v_t$ are randomly generated as $v_t = \text{normal}(\mu_t,\sigma_t)$, where $\text{normal}(a,b)$ represents a sample from a normal distribution with a mean equal to $a$ and a standard deviation equal to $b$. The temporal  factor $\gamma_t$ is computed as $\gamma_t = \prod_{\tau=1}^{t}v_{\tau}$, which is normalized as $\gamma^{\rm N}_t = \gamma_t / \max\{\gamma_t\}_{t\in\mathcal{T}}$ to satisfy that $\max\{\gamma^{\rm N}_t\}_{t\in\mathcal{T}}=1$.
    \item[4)] The load demand at each bus $b$ and time period $t$ is computed as $D_{bt}=\overline{D} \times \beta^{\rm N}_b \times \gamma^{\rm N}_t$.
    \end{itemize}

For each load profile, we solve the specific unit commitment formulation provided in the Appendix of \cite{xavier2021learning}. For the sake of simplicity, we consider neither reserves nor security constraints. Due to its high computational burden, all UCP instances are solved using the constraint generation approach proposed in \cite{fu2013modeling} using Gurobi 9.1.2 \cite{gurobi} with a MIP gap set to 0.01\%. Afterwards, the learning method is used for each instance and test system assuming that the number of neighbors is set to 50, that is, 10\% of the total number of instances. A summary of the most relevant results is provided in Table \ref{tab:results}. More specifically, columns 5 and 6 include the average and maximum suboptimality gap for each test system. Columns 7-11 contain the number of instances whose suboptimality gap belongs to given intervals. For instance, column 9 provides the number of instances with a suboptimality gap of between 0.02\% and 0.05\%. The number of instances for which all OPF problems are infeasible is included in column 12. Finally, column 13 provides the average speedup of the naive learning method in relation to solving the original UCP instances using constraint generation.

The results in Table \ref{tab:results} lead to several interesting observations. Firstly, the average suboptimality error for the four systems with more than 6000 buses is negative. This means that, on average, fixing the binary variables to those of the neighbors yields an objective function that is actually lower than that of the original UCP. Indeed, the suboptimality error is below the MIP gap in more than 99\% of the instances. Finally, these systems do not include any infeasible instances and report speedups between 41x and 172x. According to these results, we can conclude that using complicated learning-based techniques for these four systems seems unnecessary since naive alternatives significantly reduce the UCP computational burden with negligible suboptimality errors. Secondly, for systems 1888rte and 2848rte, the 50 OPF solved become infeasible for 1 and 2 instances, respectively. For these two systems, the average suboptimality errors are  slightly higher than the predefined MIP gap, while the computational times are reduced by two orders of magnitude. Therefore, the naive learning method is a competitive approach for these two systems. Finally, the number of infeasible instances for systems 1951rte, 3012wp and 3375wp are 8, 5 and 13, respectively. The average optimality errors exceed the MIP gap and amount to 0.0382\%, 0.0485\%, and 0.1256\%, in the same order, whereas the speedup factor ranges between 150x and 216x for these three systems. Consequently, using complicated learning approaches for these three systems would  only be justified when the computational time, the suboptimality error or the infeasible cases are drastically lower than those reported by the naive learning method described in this letter.

\section{Conclusion}

A highly relevant topic within the PES scientific community is the reduction of the computational burden of the unit commitment problem. However, assessing the improvements of sophisticated learning-based approaches is sometimes tricky due to the lack of benchmark methods. This letter describes a painless and naive learning method and shows that, for some power systems, learning the solution of the unit commitment problem is indeed a low-hanging fruit that can be effortlessly picked. For other systems, the results of the naive learning approach are less impressive but, at the very least, they should be used to benchmark the performance of existing and upcoming learning methods to solve the unit commitment problem.

\ifCLASSOPTIONcaptionsoff
  \newpage
\fi



\bibliographystyle{IEEEtran}
\bibliography{references}
\end{document}